\documentclass{amsart}

\usepackage[pdftex]{graphicx}
\def\l{\left}
\def\r{\right}

\def\Q{\mathbb{Q}}
\def\Z{\mathbb{Z}}

\def\e{\epsilon}
\def\C{\mathbb{C}}
\def\R{\mathbb{R}}
\def\c{\mathcal{C}}
\def\p{\mathcal{P}}
\def\dt{\Delta(\tau)}
\def\al{\alpha}

\begin{document}
\title{Integral Points on Elliptic Curves and the Bombieri-Pila Bounds}
\author{Dave Mendes da Costa}
\address{School of Mathematics, University of Bristol, University Walk, Bristol BS8 1TW, United Kingdom}
\email{madjmdc@bris.ac.uk}
\maketitle
  
\begin{abstract}
Let $\mathcal{C}$ be an affine, plane, algebraic curve of degree $d$ with integer coefficients. In 1989, Bombieri and Pila showed that if one takes a box with sides of length $N$ then $\mathcal{C}$ can obtain no more than $O_{d,\e}\left(N^{1/d+\e}\right)$ integer points within the box. Importantly, the implied constant makes no reference to the coefficients of the curve. Examples of certain rational curves show that this bound is tight but it has long been thought that when restricted to non-rational curves an improvement should be possible whilst maintaining the uniformity of the bound. In this paper we consider this problem restricted to elliptic curves and show that for a large family of these curves the Bombieri-Pila bounds can be improved. The techniques involved include repulsion of integer points, the theory of heights and the large sieve. As an application we prove a uniform bound for the number of rational points of bounded height on a general del Pezzo surface of degree 1.
\end{abstract}

\section{Background and Motivation}

Let $C$ be a plane curve (projective or affine) defined over $\Q$. The most basic number theoretic question one can ask about $C$ is to characterise the cardinality of the sets $C(\Q)$ and $C(\Z)$ (the latter being of interest when the curve is affine). An aim of Diophantine geometry is to classify the nature of these sets in terms of the geometry of $C$. This aim has seen spectacular success in the $20^{th}$ century with the genus of the curve $g(C)$ emerging as a key player.\\
\\
In this paper, we shall be concerned with integer points on affine curves (in particular the usual affine patch on an elliptic curve). The main theorem in this area is Siegel's Theorem which teaches us that if $g(C) \geq 1$ then $\# C(\Z)$ is finite. (In fact it says that this is true with $\Z$ replaced by the ring of $S$-integers of any number field.) When $C$ is $\textit{rational}$ (i.e., $g(C) = 0$), examples exist where $C(\Z)$ is finite and others where it is infinite. This behaviour is well understood.\\
\\
After one has Siegel's Theorem, a natural question is to ask for a more effective bound. Such results in this area tend to depend on the curve in question, indeed they must since there is no uniform bound for the number of integral points on a non-rational affine curve. Since one must make reference to the curve in such bounds the next step is to decide which features of the curve we wish to permit ourselves to use. This is a matter of taste and different types of bounds have different applications.\\
\\
In this paper we shall be following the tradition of trying to maintain as much uniformity as possible. As we have noted, if we desire results which are very uniform but yet still apply to a large number of curves then it is difficult to study the cardinality of $C(\Z)$. Thus attention shifts to questions concerning the $\textit{density}$ of integral points. In particular one may ask how the counting function

\[ C(\Z,N) := \# \{(x,y) \in C(\Z) : \left|x\right|, \l|y\r| \leq N \} \]
\\
behaves as $N>0$ gets very large. Since we want our results to be uniform, an asymptotic result is not the correct thing to be looking for. Instead we would like a bound on this function which depends on $N$ and on as little to do with the curve as possible. If we were to look to bound $C(\Z,N)$ independently of any information about $C$ then the best we could do would be the trivial bound $C(\Z,N) \ll N^2$ since there are curves of high degree which subsume all the integer points in the box. Therefore we should ask for a bound in terms of the degree of the curve.\\
\\
The breakthrough result in this area was reported in 1989 when Bombieri and Pila \cite{BP} showed that if $C$ is a plane curve of degree $d$ then 

\begin{equation}\label{eq:1} C(\Z,N) \ll_{d,\e} N^{1/d + \e} .\end{equation}
\\
This bound is tight as is seen by considering the curves

\[ C_d : y = x^d \] 
\\
The interesting thing about this family is that it consists of rational curves only. Indeed, the only known examples where ~\eqref{eq:1} is tight is when $C$ is rational. From the point of view of Diophantine geometry this is intriguing since it once again points to differing arithmetic behaviours being catagorised by geometry. Indeed, this prompts us to ask the following question:

\newtheorem{Question}{Question}
\begin{Question}
If $C$ is a non-rational curve of degree $d$ then is there a $\delta(d) >0$ depending only upon $d$ such that

\[ C(\Z,N) \ll_d N^{1/d - \delta(d)} ? \]
\end{Question}

Currently, the best evidence that this question has a positive answer is that the corresponding question for rational points does. In 2000, Heath-Brown \cite{HB} showed that if one takes a plane projective curve of degree $d$ and counts points of (naive) height less than $N$ (that is, points $P = [x:y:z]$ such that $x,y,z$ are coprime integers of absolute value less than or equal to $N$) then for the associated counting function one gets

\[ C(N,\Q) \ll_{d,\e} N^{2/d + \e} .\]

Once again the only examples where this is tight are the projective versions of the curves we mentioned before and so the same question as above can be asked about non-rational curves. This was answered in the positive in 2005 by Ellenberg and Venkatesh \cite{EV}. So there is hope that our question has a positive answer and this answer is in reach even though the methods in \cite{EV} do not directly apply.\\
\\
This paper is concerned with the case where $C$ is a non-singular degree 3 non-rational curve, i.e., an elliptic curve. We shall show that for a large class of such curves we can give a positive answer to the Question. In particular we show the following two results. Firstly we have our

\newtheorem*{main}{Main Theorem}
\begin{main}
\label{th:main}
Let $E$ be an elliptic curve given by the equation $E: y^2 = f(x)$ where $f \in \Z[x]$ is a monic cubic with no repeated roots. Then there is a $\delta > 0$, independent of $E$, such that $E(\Z,N) \ll N^{1/3-\delta}$ (the implied constant is uniform).\\
\end{main}

Our second result strengthens and extends this to a wider class of elliptic curves but requires more hypotheses:

\newtheorem{thm}{Theorem}
\label{th:thm1} 
\begin{thm}
Let $E$ be an elliptic curve given by the equation $y^2 = f(x)$ where $f \in \Z[x]$ is a monic cubic with no repeated roots. Suppose further $c_4(E) < 0$ and $j(E) > \e$ for some $\e > 0$. Let $\mathcal{B}$ be any box in the plane with sides of length $N$. Then there is a $\delta > 0$ (depending only on $\e$) such that $\# \left(E(\Z) \cap \mathcal{B}\right) = O_{\e}(N^{1/3-\delta})$.  
\end{thm}

\subsection{Acknowledgments} I am glad to acknowledge the useful conversations and insights that many people have shared with me. In particular, I extend my warmest regards to Tim Browning, Tim Dokchitser, Roger Heath-Brown, Marc Hindry, Pierre Le Boudec, Dan Loughran, Joe Silverman and Trevor Wooley. Many thanks go to my advisor Harald Helfgott who gave me this problem and has always been generous both as a supervisor and a friend. I would also like to thank the University of Bristol and the Ecole Normale Superieure for providing me with wonderful working conditions during my PhD. Finally, dues are paid to the EPSRC Doctoral Training Award, which has funded my research, with the hope that they may continue to support such pure mathematical endevours. 

\section{Notations and Conventions}

We shall let $E$ denote an elliptic curve throughout. Our problem is phrased with respect to particular embeddings of $E$ and so we shall be assuming that the model of $E$ takes a Weierstrass form. We say that $E$ is in $\textit{long Weierstrass form}$ when we have

\[ E: y^2 + a_1 xy + a_3 y = x^3 + a_2 x^2 +a_4 x + a_6 \]
\\
and in $\textit{short Weierstrass form}$ if we have

\[ E: y^2 = x^3 + Ax + B \] 
\\
where $a_i,A,B \in \Z$. We make no assumptions on the minimality of $E$. We also associate several standard quantities with $E$, namely the discriminant $\Delta_E$, $c_4(E)$ and the $j$-invariant $j(E)$. In the case of a curve in short Weierstrass form these are:

\[ \Delta_E = -16(4A^3+27B^2), \text{  } c_4(E) = -\frac{A}{27}, \text{  } j(E) = -1728\frac{(4A)^3}{\Delta_E} .\]
\\
In our bounds we shall make use of both Landau's Big-Oh notation and Vinogradov's $\ll$ notation. Given functions $f,g: \R \rightarrow \R$ we shall write either 

\[ f(x) = O(g(x)) \text{  or  } f(x) \ll g(x) \]
\\
in the case that there is a constant $C > 0$ such that $\left|f(x)\right| \leq C \left|g(x)\right|$. In particular $f = O(1)$ implies $f$ is bounded. The constant $C$ is suppressed by the notation and so is referred to as the $\textit{implied constant}$. If we wish to draw note to the dependency of the implied constant on something (such as the curve $E$) then we shall do so in subscript writing, for instance,

\[ f(x) = O_E (g(x)) .\]
\\
We shall make reference to several different heights in what follows. Let $P = (x,y) \in \Z^2$. By the $\textit{naive height}$ of $P$ we shall mean the maximum of $\l|x\r|$ and $\l|y\r|$. When dealing with points on an elliptic curve we shall be studying the canonical height $\hat{h}$ on $E$. Despite its canonical nature, there are still several ways of scaling this height which are used in the literature. We shall be using the scaling which corresponds to seeing $\hat{h}$ as a sum of local heights as defined in \cite{sil}. This is half the size of the version defined by Tate.

\section{Repulsion techniques}

The majority of the work in this paper goes into putting us into a situation where we can use a result about how integral points can be said to $\textit{repel}$ one another. The idea of rational points repelling each other has a wonderful pedigree. In 1969, Mumford showed that if $C$ is a curve of genus $g>1$ then two rational points $P, Q \in C(\Q)$, which in the Jacobian of $C$ have comparable canonical height, repel each other in the sense that the angle between them in the Mordell-Weil lattice of Jac$(C)$ is at least $\arccos(\frac{1}{g})$. This observation, known as Mumford's Gap Principle, allowed Mumford to show that the number of rational points on $C$ of naive height less than $N$ is 

\[ O_C(\log \log N) .\]
\\
This result was improved by Faltings who showed that 

\[ \# C(\Q) = O_C(1) .\]
\\
This bound can also be seen as a consequence of repulsion, the key difference being that the condition that the points have to be of comparable canonical height is shown to not be necessary to force repulsion. So repulsion has a lot to say about rational points; how about integral points?\\
\\
In the case of an elliptic curve, Mumford's Gap Principle says nothing about rational points other than they repel each other by an angle of at least 0! However, if we restrict to integral points then we obtain a non-trivial repulsion of at least $60^\circ$ between points which have comparable canonical height. This phenomenon is present in the work of Silverman who showed how to integrate the notion of (S-)integrality of points into such repulsion techniques in \cite{sil2} (later refined with Gross in \cite{GS}). If we demand that the points reduce to the same point modulo some prime $p$, we can increase the minimum angular repulsion between points at the cost of having to partition the points into the fibres of the reduction map. To turn this repulsion into a bound on the number of integral points of canonical height less than $h$, we simply slice up the set of integral points into strips with comparable canonical height, i.e.,

\[ (1-\epsilon)\log N \leq \hat{h}(P) \leq (1+\e)\log N \]
\\
for some $N>0$, and then to each slice apply bounds for the number of vectors which can packed with suitable angular repulsion onto a sphere of dimension $r - 1$ where $r$ is the rank of $E(\Q)$. The idea of gaining bounds through slicing and sphere-packing was first seen in the thesis of Helfgott and in the paper \cite{Helf}. These ideas culminated in the paper \cite{HV} of Helfgott and Venkatesh who show how this repulsion can be obtained in a uniform manner and also how the choice of the prime $p$ can be optimised to yield the following bound:

\begin{thm}\cite[Corollary~3.9]{HV}
\label{th:thm2} 
Let $E$ be an elliptic curve over $\Q$ with integer coefficients and $\e > 0$. Then there is a $\delta > 0$ such that the number of integral points of $E$ of canonical height less than $h$ is 
\[ O_\e(e^{h(1-\delta)+\e}) \]
\\
where the implied constant does not depend on $E$.\\
\end{thm}

Thus we see that this theorem yields the desired result when $h \leq \frac{1}{3}\log N$ and in fact we can let $h$ be a touch larger. In particular, the key to our application of Theorem ~\ref{th:thm2} will be the following result.

\newtheorem{Prop}{Proposition}

\begin{Prop}
\label{th:thm3}
Let $E$ be an elliptic curve in short Weierstrass form and $N,\delta>0$ such that $\left|\Delta_E\right| < N^{4+6\delta}$ and $\left|A\right|\leq N^{4/3+2\delta}$. Then if $Q = (x,y) \in E(\Z)$ has $\left|x\right| \leq N^{2/3+\delta}$ then 
\[ \hat{h}(Q) \leq \left(\frac{1}{3}+\frac{\delta}{2}\right)\log N + O(1) \] where the implied constant is absolute. 
\end{Prop} 

\section{Proof of Proposition ~\ref{th:thm3}}

The idea of the proof is to break up the canonical height, $\hat{h}$, as a sum of local canonical heights $\lambda_p$, one for each prime $p$ and then bound each part. We shall use the local heights as they are defined in \cite[Chapter~VI]{sil}.\\
\\
For the finite primes we shall use the following bound, due to Tate.

\newtheorem{lem}{Lemma}
\begin{lem}
\label{th:lem1}
Let $p$ be a finite prime, $\lambda_p$ the local height of $E$ at $p$ and let $Q \in E(\mathbb{\Z})$. Then
\[ \lambda_p(Q) \leq \frac{\text{ord}_p(\Delta_E)}{12}\log p .\]
\end{lem}

\proof{\cite[Chapter~II, Theorem 4.5]{La}}\\
\qed\\
\\
Thus we see that for $Q \in E(\Z)$,

\begin{eqnarray}\label{eq:tate} \sum_{p \neq \infty}{\lambda_p(Q) \leq \frac{1}{12}\log \left|\Delta_E \right|} .\end{eqnarray}
\\
This leaves just the infinite part to estimate. Let us refresh ourselves on the theory of elliptic curves over $\mathbb{C}$. \\
\\
Every elliptic curve $E$ over $\mathbb{C}$ is isomorphic to $\mathbb{C} / \Lambda$ for some rank 2 lattice $\Lambda \subset \C$. The isomorphism is explicit and comes from the Weierstrass $\wp_\Lambda$-function which maps $\C / \Lambda$ to $E(\C)$ and is defined by

\[ \wp_\Lambda(z) = \frac{1}{z^2} + \sum_{k>0}{(2k+1)G_{2k+2}(\Lambda)z^{2k}}\]
\\
where 

\[G_{2k}(\Lambda) = \sum_{0\neq \omega \in \Lambda}{\frac{1}{\omega^{2k}}} \]
\\
is the weight $2k$ Eisenstein series associated to $\Lambda$. The isomorphism between $\C / \Lambda$ and $E(\C)$ is given by

\[ z \mapsto (\wp_\Lambda(z),\wp^\prime_\Lambda(z)/2) .\]
\\
In this case the equation for $E$ can be recovered from $\Lambda$ via 

\[ E: y^2 = x^3 - 30G_4(\Lambda)x - 70G_6(\Lambda) .\]
\\
Two curves $\C / \Lambda_1$ and $\C / \Lambda_2$ are isomorphic over $\C$ if and only if there is some $u \in \C^\times$ such that $u\Lambda_1 = \Lambda_2$. In this way one can normalise the form of $\Lambda$ to be $\Lambda_\tau = <1,\tau>$ where 

\[ \tau \in \left\{z \in \C : \left|\text{Re}(z)\right| \leq \frac{1}{2} , \left|z\right| \geq 1 \right\} =: \mathcal{F} .\]
\\
For such a $\tau$ let $E_\tau$ be the curve associated to $\Lambda_\tau$.\\
\\
Another way of characterising two curves being isomorphic is via the $j$-invariant 

\[ j_E = -1728 \frac{(4A)^3}{\Delta_E} .\]
\\
Two curves, $E_1$ and $E_2$ are isomorphic if and only if $j_E = j_{E^\prime}$. Moreover, if $E$ is defined over a field $K$ then $j_E \in K$. The curves we are interested in are defined over $\Q$ and thus over $\R$, so it is sensible to ask which curves in $\mathcal{F}$ have real $j$-invariant. 

\begin{lem}
Let $E/\R$ be an elliptic curve. Then there is a unique $\tau$ in the set 

\[ \mathcal{C} = \left\{ ib : b \geq 1 \right\}\cup \left\{e^{i\theta} : \frac{\pi}{3} \leq \theta < \frac{\pi}{2}\right\} \cup \left\{\frac{1}{2} + ib : b > \frac{\sqrt{3}}{2} \right\}  \]
\\
such that $j(E_\tau) = j(E)$.
\end{lem}
\proof{\cite[Chapter~V, Prop~2.1, p.414]{sil}} \qed\\
\\
If $j(E) = j(E_\tau)$ then we shall say that $E$ is $\textit{associated to } \tau$ and if the isomorphism between $E$ and $E_\tau$ is defined over a field $K$ we shall say it is $\textit{associated to } \tau \textit{ over }$ $K$. Let us label the three sets which make up $\c$ by $\c_1, \c_2$ and $\c_3$ respectively. These regions neatly separate up the elliptic curves defined over $\R$. The curves which are associated to $\tau \in \c_1$ are those with two real components. These curves have $A \leq 0$ and $j \geq 1728$. The curves associated to $\tau \in\c_3$ have one real component, $A < 0$ and $j < 0$. Finally, the curves associated to $\tau \in \c_2$ have one real component, $A > 0$ and $0 \leq j < 1728$. It will be noted that no mention of the $B$ coefficient has been made. This is because we can change the sign of the $B$ coefficient via the isomorphism

\[ (x,y) \mapsto (-x,iy) \]
\\
It is worth noting that in terms of the underlying lattices, this isomorphism is 

\[ \Lambda \mapsto \frac{1}{i}\Lambda \]
\\
and so consists of rotating the lattice through $90^\circ$. We shall refer to the image of $E_\tau$ under this map as its $\textit{twist}$. We note that every elliptic curve defined over $\mathbb{R}$ is isomorphic over $\mathbb{R}$ to either $E_\tau$ or its twist for some $\tau \in \mathcal{C}$. \\
\\
We shall now look at the set $E(\R)$ in a bit more detail. Let $E(\C) \cong E_\tau(\C)$ and consider a fundamental parallelogram for $\Lambda_\tau$

\[ \p_\tau = \left\{ x_1 + x_2\tau : x_1,x_2 \in (-1/2,1/2] \right\} .\]
\\
For every elliptic curve $E$ which has $j(E) \in \R$, there is a $u \in \C^\times$ such that $\overline{u\Lambda} = u\Lambda$. Let $\Lambda^\prime = u\Lambda$. Then $\overline{\wp_{\Lambda^\prime}(z)}=\wp_{\Lambda^\prime}(\bar{z})$ and so the real points are those which are invariant under the action of complex conjugation on $\C / \Lambda^\prime$. \\
\\
For curves $E_\tau$ with $\tau \in \c_1$ or $\c_3$ we already have $\overline{\Lambda_\tau}=\Lambda_\tau$. This means that in the case of $\c_1$, the (two) real components are the images of the points 

\[ \left\{ x_1 + x_2\tau \in \p_\tau : x_2 \in \left\{0,\frac{1}{2}\right\} \right\} \]
\\
and for $\c_3$ the (one) real component is the image of

\[ \left\{ x_1 + x_2\tau \in \p_\tau : x_2 = 0 \right\} .\]
\\
In the case of $\c_2$ we need to modify $\Lambda_\tau$ to $\Lambda^\prime = \frac{1}{\tau^{1/2}}\Lambda$. This has the effect of rotating the point $1 + \tau$ onto the real axis. Thus the real component is the image of 

\[ \left\{ t(1 + \tau) : -\frac{1}{2} < t \leq \frac{1}{2} \right\} .\]
\\
The final thing left to note is that if we take the twist of $E_\tau$ (which flips the sign of the $B$ coefficient), call this twist $E_\tau^\prime$, then since this rotates the lattice by $90^\circ$ it is easy to find the points of $\p_\tau$ which correspond to points on $E^\prime (\R)$ (after we rotate the lattice). They are the points which, after we rotate $\Lambda_\tau^\prime$ by $90^{\circ}$, map to themselves under complex conjugation. These are the points on $E(\mathbb{C})$ which are mapped to $E^\prime(\R)$ by the twist.\\
\\
Let $\tau \in \c$ and let $E$ be associated to $\tau$ over $\R$. We want to study $\lambda_\infty(u)$ for those $u \in \p_\tau$ which are mapped to $E(\R)$ and where the naive height of the $x$-coordinate is bounded. We have now identified the points which map to $E(\R)$ and later we shall find conditions on $u$ which correspond to the bound on the $x$-coordinate. We also wish to study the same points for the twist $E^\prime$ of $E$. Rather than studying $\lambda_\infty$ on the fundamental domain for $E^\prime$, it is sufficient for us to understand the behaviour of $\lambda_\infty(u)$ for those $u \in \p_\tau$ which correspond to points on $E^\prime (\R)$ after we twist. This is sufficient since the value of $\lambda_\infty(u)$ does not change as we twist $E$ to $E^\prime$. Therefore we shall, in the sequel, study the behaviour of $\lambda_\infty(u)$ for those $u \in \p_\tau$ which map to either $E(\R)$ or to the image of $E^\prime(\R)$ in $E(\mathbb{C})$.\\     
\\
We now have enough background to start examining the archimedean local canonical height.\\
\\
Take a point $u = u_1 + u_2\tau \in \p_\tau$, then following \cite[p.~468]{sil} we can define the archimedean local canonical height to be

\[ \lambda_\infty(u) = -\frac{1}{2}B_2(u_2)\log\left|q\right| -\log\left|1-t\right| - \log\left|\prod_{n>0}(1-q^nt)(1-q^nt^{-1})\right| \]
\\
where $q = e^{2\pi i \tau}$, $B_2(x) = x^2 -x + 1/6$ is the second Bernoulli polynomial and $t = e^{2 \pi i u}$. \\
\\
In order to analyse the behaviour of $\lambda_\infty$ for points of $E^{(\prime)}(\R)$, we note that there is a similarity between certain terms in $\lambda_\infty$ and the Jacobi product formula for the discriminant function when viewed as a function on $\mathcal{F}$,

\begin{equation}
\label{jacobi} \Delta(\tau) = \Delta(E_\tau) = (2\pi)^{12}q\prod_{n>0}(1-q^n)^{24} .\end{equation}
\\
We can use this similarity to our advantage.

\begin{lem}
\label{th:lem3}
Let $\tau \in \mathcal{F}$. Then 
\[ \lambda_\infty(u) \leq -\log\left|1-t\right| - \frac{1}{12}\log\left|\Delta(\tau)\right| + O(1) \]
\\
where the implied constant in the $O(1)$ is absolute.
\end{lem}
\proof{First we note that since $\lambda_\infty(u) = \lambda_\infty(-u)$ for all $u \in \p_\tau$ we can choose $u$ such that $u_2 \in [0,1/2]$. Noting that $\left|q\right| < 1$ and that $B_2(u_2) \in [-1/12, 1/6]$ we have

\[ -\frac{1}{2}B_2(u_2)\log\left|q\right| \leq -\frac{1}{12}\log\left|q\right| .\]
\\
This corresponds to the $q$ term in $ \Delta(\tau)$. Next we note that $t = e^{2\pi i u} = e^{2 \pi i u_1}q^{u_2}$, and so

\begin{eqnarray*}
\left|\prod_{n>0}(1-q^nt)(1-q^nt^{-1})\right| & = & \left|\prod_{n>0}(1-e^{2 \pi i u_1} q^{n+u_2})(1-e^{-2 \pi i u_1} q^{n-u_2})\right| \\
								                 & \geq & \prod\limits_{n>0} (1-\left|q\right|^{n+u_2})(1-\left|q\right|^{n-u_2}) \\
								                 & \geq & \prod\limits_{n>0} (1-\left|q\right|^{n})(1-\left|q\right|^{n-1/2}) \\
								                 & \geq & \prod\limits_{n>0}(1-e^{-\sqrt{3}\pi n})(1-e^{-\sqrt{3}\pi(n-1/2)}) \\
								                 & = & 0.92984\ldots 
\end{eqnarray*}
\\
where in the penultimate line we use the fact that $\text{Im}(\tau) \geq \frac{\sqrt{3}}{2}$. Therefore $\lambda_\infty(u) \leq -\frac{1}{12}\log\left|q\right| - \log\left|1-t\right| + O(1)$, which only leaves us to show that 

\[ \log\left|\Delta(\tau)\right| = \log\left|q\right| + O(1) \]
\\
and this is the case since the lower bound on $\text{Im}(\tau)$ allows us to deduce that

\[ 21.588 \leq \log\left|(2\pi)^{12}\prod_{n>0}(1-q^n)^2\right| \leq 22.4554 .\]
\\
Thus we have our result. 

\qed \\}
\\
Here is the strategy for the rest of the proof. Let $E$ be given in short Weierstrass form and suppose that $\left|\Delta_E\right| \leq N^{4+6\delta}$ and $\left|A\right	|\leq N^{4/3+2\delta}$. There is a $\tau \in \mathcal{C}$ such that $E$ and $E_\tau$ are isomorphic over $\mathbb{C}$. We want to understand the behaviour of $\lambda_\infty$ at points of $\p_\tau$ which correspond to points $Q = (x(Q),y(Q)) \in E(\R)$ and satisfy $\left|x(Q)\right| \leq N^{2/3+\delta}$. Since $E$ and $E_\tau$ are in short Weierstrass form we know that there is a $w \in \mathbb{C}$ such that

\begin{eqnarray}\label{eq:2}
E_\tau(\mathbb{C}) & \stackrel{\cong}{\rightarrow} & E(\C) \nonumber\\
(x,y) & \mapsto & (w^2x, w^3 y) 
\end{eqnarray}
\\
and so we have $x(Q) = w^2\wp_\tau(u)$ since $(\wp(u),\wp^\prime (u))$ parameterises $E_\tau(\mathbb{C})$. The idea of the rest of the proof is as follows: we can deduce the size of $w$ from the ratio of $\Delta_E$ and $\Delta(\tau)$ since the discriminant is a weight 12 modular form and so the ratio is $w^{12}$. Next, we use this to get an upper bound on the size of $\left|\wp(u)\right|$ which is necessary for $\left|x(Q)\right|$ to be smaller than $N^{2/3+\delta}$. Finally we can use this upper bound to get a lower bound on $\left|u\right|$.  This bound will turn out to be just what we want in order to cancel out the contribution from the local heights at the finite primes.\\
\\
We shall start by understanding $\wp_\tau$ a bit better and since

\[ \wp_\tau(z) = \frac{1}{z^2} + \sum_{k\geq 1}{(2k+1)G_{2k+2}(\tau)z^{2k}} \]
\\
it makes sense to start with the $G_{2k}$. 

\begin{lem}
\label{th:lem4}
For $\tau \in \c$ and $k\geq 2$ we have $\left|G_{2k}(\tau)\right| \leq 80$.
\end{lem}                                                                                                  

\proof{We have

\[\left|G_{2k}(\tau)\right| = \left|\sum_{0\neq\omega \in \Lambda_\tau}{\frac{1}{\omega^{2k}}}\right| 
								                  \leq  \sum{\left|\frac{1}{\omega^{2k}}\right|} 
								                  \leq  \sum{\left|\frac{1}{\omega^{4}}\right|} \]
\\								                  
where the last line is justified by the fact that, for $\tau$ in $\c$, the $\left|\omega\right| \geq 1$ for every $\omega \in \Lambda_\tau$. Thus it only remains to bound this last term for all $\tau \in \c$.\\
\\
For $\tau = bi \in \c_1$ we have
\begin{eqnarray*}
\left|G_{2k}(\tau)\right| & \leq & \sum_{0\neq\omega \in \Lambda_\tau}{\left|\frac{1}{\omega^{4}}\right|} \\
                          & = & \sum_{(m,n)\in\Z^2-\{0\}}{\frac{1}{\left|m+nbi\right|^4}} \\
                          & = & \sum{\frac{1}{(m^2+(nb)^2)^2}} \\
                          & \leq & \sum{\frac{1}{(m^2+n^2)^2}} \\
												  & \leq & 7 .
\end{eqnarray*}
Similarly, if we consider $\tau = \frac{1}{2} + bi$ with $b \geq \frac{1}{2}$ then we have

\begin{eqnarray*}
\left|G_{2k}(\tau)\right| & \leq & \sum_{0\neq\omega \in \Lambda_\tau}{\left|\frac{1}{\omega^{4}}\right|} \\
                          & = & \sum_{(m,n)\in\Z^2-\{0\}}{\frac{1}{\left|m+(\frac{1}{2}+bi)n\right|^4}} \\
                          & = & \sum{\frac{1}{(m^2+((\frac{1}{2}+bi)n)^2)^2}} \\
                          & = & \sum{\frac{1}{(m^2+(nb)^2 + n^2/4)^2}} \\
												  & \leq & \sum{\frac{1}{(m^2+(\frac{n}{2})^2)^2}} \\
												  & \leq & 80. 
\end{eqnarray*} 
\\
This just leaves $\tau \in \c_2$. For this we note that the transformation

\[ \alpha : \tau \mapsto \frac{-1}{\tau - 1} \]
\\
takes $\c_2$ to the line segment

\[ \left\{ \frac{1}{2}+bi : \frac{1}{2}\leq b\leq \frac{\sqrt{3}}{2} \right\} .\]
\\
Since $G_{2k}$ is a weight $2k$ modular form we have the relation

\begin{eqnarray*}
\left|G_{2k}(\tau)\right| & = & \left|\tau - 1\right|^{-2k}\left|G_{2k}(\alpha(\tau))\right| \\
                           & \leq & \left|G_{2k}(\alpha(\tau))\right| \\
                           & \leq & 80.
\end{eqnarray*}

Thus $\left|G_{2k}(\tau)\right| \leq 80$ for all $\tau \in \c$.

\qed
\\
With this estimate we learn the following fact about $\wp_\tau$:

\newtheorem{cor}{Corollary}
\begin{cor}
\label{th:cor1}
Let $\tau \in \c$ and $u \in \p_\tau$. Then 
\[ \left|\wp_\tau(z)\right| \geq \frac{1}{\left|z\right|^2} - 100.\]
\end{cor}
\proof{First note that from the previous lemma and the fact that $\left|u\right| \leq \frac{1}{2}$, since $u \in \mathcal{P}_\tau$, we have
\begin{eqnarray*}
\left|\sum_{k>1}{(2k+1)G_{2k+2}(\tau)u^{2k}}\right| & \leq & \sum_{k>1}{(2k+1)\left|G_{2k+2}(\tau)\right|\left|u\right|^{2k}} \\
                                                    & \leq & 80\sum_{k>1}{(2k+1)\left|u\right|^{2k}} \\
                                                    & \leq & 80\sum_{k>1}{(2k+1)(0.5)^{2k}} \\
                                                    & = & 97.7778 \\
                                                    & \leq & 100 
\end{eqnarray*}
\\
and so we have

\begin{eqnarray*}
\left|\wp_\tau(z)\right| & \geq & \frac{1}{\left|z\right|^2} - \left|\sum_{k>1}{(2k+1)G_{2k+2}(\tau)z^{2k}}\right| \\
                         & \geq & \frac{1}{\left|z\right|^2} - 100. 
\end{eqnarray*} 
(Of course this bound is trivial unless $\left|u\right| \leq \frac{1}{\sqrt{100}}$.)
\qed }\\
\\
We shall now restrict our attention to points $u \in \p_\tau$ which map to $E_\tau(\R)$ or $E^{\prime}_\tau(\R)$ (after the twist). Our aim is to use the above bound to get a lower bound on $\left|u\right|$ for those $u \in \p_\tau$ which map to points $Q$ on $E(\R)$ or $E^\prime(\R)$ such that $\left|x(Q)\right| \leq N^{2/3+\delta}$. Since the result of Corollary ~\ref{th:cor1} is trivial unless $\left|u\right| \leq \frac{1}{\sqrt{100}}$ we shall need to make a hypothesis giving an upper bound on $\left|u\right|$. 

\begin{lem}
\label{th:ubound}
Consider $\tau \in \c$ and let $E$ be the elliptic curve associated to $\tau$ over $\R$. Consider those $u \in \p_\tau$ which map to points $Q = (x(Q),y(Q))$ in $E(\mathbb{C})$ such that $\left|x(Q)\right| \leq N^{2/3+\delta}$. Then if $\left|u\right| \leq \frac{1}{\sqrt{200}}$ we have that
\[ \left|u\right| \geq \frac{1}{\sqrt{2}}N^{\eta/12 - 1/3-\delta/2}\left|\Delta(\tau)\r|^{-1/12} \]
\\
where $\eta > 0$ is such that $\left|\Delta_E\right| = N^{\eta}$. 
\end{lem}

\proof{
We shall argue by the contrapositive. Thus we wish to find a value $\beta(N) > 0$ such that if $\left|u\right| < \beta(N)$ then we have $\left|x(Q)\right| > N^{2/3+\delta}$. Since $E$ and $E_\tau$ are in short Weierstrass form we know that there is a $w \in \mathbb{C}$ such that

\[ \left|w^{12} \Delta(\tau)\right| = \left|\Delta_E\right| = N^\eta. \]
\\
Thus

\[ \left|w^2\right| = N^{\eta/6}\left|\dt\right|^{-1/6} .\]
\\
This means that

\begin{eqnarray*}
 \left|x(Q)\right| & = & \left|w^2\wp(u)\right| \\
                   & = & N^{\eta/6}\left|\dt\right|^{-1/6}\left|\wp(u)\right|
\end{eqnarray*}
\\
and so if we want $\left|x(Q)\right| > N^{2/3+\delta}$ then 

\begin{eqnarray}\label{eq:3}
N^{2/3+\delta} & < & N^{\eta/6}\left|\dt\right|^{-1/6}\left|\wp(u)\right| 
\end{eqnarray}
\\
is sufficient. Now, since $\left|u\right| \leq \frac{1}{\sqrt{200}}$ then Corollary ~\ref{th:cor1} tells us that

\[ \left|\wp(u)\right| \geq \frac{1}{\left|u\right|^2} - 100 \geq \frac{1}{2\left|u\right|^2} .\]
\\
To ensure that condition ~\eqref{eq:3} is met we have the stronger condition:

\[ N^{2/3+\delta} < \frac{1}{2}N^{\eta/6}\left|\dt\right|^{-1/6}\frac{1}{\left|u\right|^2} \]\\
and so we see that $\left|x(Q)\right| > N^{2/3+\delta}$ if

\[ \left|u\right| < \frac{1}{\sqrt{2}}N^{\eta/12-1/3-\delta/2}\left|\dt\right|^{-1/12} \].
\qed}
\\
The final step is to find an upper bound on $-\log\left|1-t\right|$. We start by making a small but useful observation.

\begin{lem}
\label{th:lem5}
Let $x \in [0,1]$, then 

\[ \left|1-e^{ix}\right| \geq \frac{x}{2} \]
\\
and 
\[ \left|1-e^{-x}\right| \geq \frac{x}{2} .\]
\end{lem}  
\proof{This is easily seen after expanding the power series for $e^{-x}$ and $e^{ix}$.}
\qed
\begin{lem}
\label{th:lem6}
Let $u \in \p_\tau$ be a point which maps to $E(\R)$ or the image of $E^\prime(\R)$ in $E(\mathbb{C})$ and suppose $u$ satisfies  $\frac{1}{\sqrt{2}}N^{\eta/12 - 1/3-\delta/2}\left|\Delta(\tau)\r|^{-1/12} < \left|u\right| < \frac{1}{\sqrt{200}}$ , then 
\[ -\log|1-t| \leq \left(\frac{1}{3}+\frac{\delta}{2}-\frac{\eta}{12} \right)\log N + \frac{1}{12}\log\left|\dt\right| + O(1) \]
\\
where the $O(1)$ is absolute.
\end{lem}

\proof{Let $\tau =  bi \in \c_1$ or $\frac{1}{2} + bi \in \c_3$. Since we are interested in the points of $\p_\tau$ which map to a point in either $E(\R)$ or $E^\prime(\R)$ we recall from our previous discussion of such points that $u$ is either in the interval $[0,1/2]$ or $[0,bi/2]$ depending on the sign of $B$ in the equation for $E$. (Note that we are using the assumption that $\l|u\r| \leq \frac{1}{\sqrt{200}}$ here to ensure that we are on the connected component containing the identity element.) If we are on the interval $[0,1/2]$ then we have

\begin{eqnarray*}
-\log\left|1-t\right| & = & -\log\left|1-e^{2\pi i u}\right| \\
                      & \leq & -\log\left|\pi u\right| \\
                      & \leq & -\log\left|\pi N^{\eta/12-1/3-\delta/2}\dt^{-1/12}\right|\\
                      & = & \left(\frac{1}{3}+\frac{\delta}{2}-\frac{\eta}{12}\right)\log N + \frac{1}{12}\log\left| \dt \right| + O(1)
\end{eqnarray*}
\\
where the second line is justified by Lemma ~\ref{th:lem5}. Next, let $u \in [0,bi/2]$. Then we have     

\begin{eqnarray*}
-\log\left|1-t\right| & = & -\log\left|1-e^{-2\pi \left|u\right|}\right| \\
                      & \leq & -\log\left|\pi u\right| \\ 
                      & \leq & -\log\left|\frac{\pi}{\sqrt{2}} N^{\eta/12-1/3-\delta/2}\dt^{-1/12}\right| \\
                      & = & \left(\frac{1}{3}+\frac{\delta}{2}-\frac{\eta}{12}\right)\log N + \frac{1}{12}\log\left| \dt \right| + O(1)
\end{eqnarray*}                  
\\
where the second line follows from the fact that $2 \pi \left|u\right| < \frac{2\pi}{\sqrt{200}} < 1$ and Lemma ~\ref{th:lem5}.\\
\\
This gives us our result in this case. Only $\c_2$ remains. Let $\tau = e^{i\theta}$ then we have $u$ on the line

\[ \left\{ x(1+\cos \theta + i \sin \theta) : -\frac{1}{2} < x \leq \frac{1}{2} \right\} .\]
\\
This means that we have

\[ \l|u\r| = \l|x\r|\sqrt{2(1+\cos \theta)}.\]
\\
We can assume that $x \geq 0$ and so we have

\begin{eqnarray*}
-\log\left|1-t\right| & = & -\log\left|1-e^{2\pi i x(1+\cos \theta + i \sin \theta)}\right| \\
                     & \leq & -\log\left|1-e^{-2\pi x \sin \theta}\right| \\
                     & = & -\log\left|1-e^{\frac{-2\pi \left|u\right| \sin \theta}{\sqrt{2(1+\cos \theta)}}}\right| .
\end{eqnarray*}                     

Since $\tau \in \mathcal{C}$ we have $\frac{\pi}{3} \leq \theta \leq \frac{\pi}{2}$ and thus $\frac{\sqrt{3}}{2} \leq \sin \theta$ and $\cos \theta \geq 0$. Hence

\begin{eqnarray*}                     
-\log\left|1-t\right|& \leq & -\log\left|1-e^{-\pi\sqrt{3}\left|u\right|}\right| \\
                     & \leq & -\log\left|\sqrt{3} \pi u\right|\\
                     & \leq & -\log\left|\frac{\sqrt{3}\pi}{\sqrt{2}}N^{\eta/12-1/3-\delta/2}\dt^{-1/12}\right|\\  
                     & \leq & \left(\frac{1}{3}+\frac{\delta}{2}-\frac{\eta}{12}\right)\log N + \frac{1}{12}\log\left| \dt \right| + O(1)
\end{eqnarray*} 

and so we have shown the lemma to be true.}
}
\qed\\
\\
\begin{cor}
\label{th:cor2}
Let $E$ be an elliptic curve over $\mathbb{Q}$ given by $E: y^2 = x^3 + Ax + B$. Let $N>0$ and $Q= (x(Q),y(Q)) \in E(\Z)$ with $\left|x(Q)\right| < N^{2/3+\delta}$. Suppose further that $E$ is associated with $\tau$ and $Q$ corresponds to $u \in \p_\tau$ with $\left|u\right| < \sqrt{\frac{1}{200}}$. Then the statement of Proposition ~\ref{th:thm3} is true in this case, i.e.,

\[ \hat{h}(Q) \leq \left(\frac{1}{3}+\frac{\delta}{2}\right)\log N + O(1) \] 
\\
where the implied constant is absolute.
\end{cor}

\proof{ Let $\eta > 0$ be such that $\left|\Delta_E\right| = N^{\eta}$. Then
\begin{eqnarray*}
\hat{h}(Q) & = & \sum_{p\leq \infty}\lambda_p(Q) \\
           & \leq & \frac{1}{12}\log \left|\Delta_E\right| + \lambda_\infty(u) \\
           & \leq & \frac{\eta}{12}\log N - \log\left|1-t\right| - \frac{1}{12}\log\left|\dt\right| + O(1)
\end{eqnarray*}
\\
where the second line follows from the inequality ~\eqref{eq:tate} derived from Lemma ~\ref{th:lem1} and last line is justified by Lemma ~\ref{th:lem3}. Now, since $\left|x(Q)\right|<N^{2/3+\delta}$ we see that Lemma ~\ref{th:ubound} and Lemma ~\ref{th:lem6} can be combined to obtain
\begin{eqnarray*}
\hat{h}(Q) & \leq & \frac{\eta}{12}\log N + \left(\frac{1}{3}+\frac{\delta}{2}-\frac{\eta}{12} \right)\log N + \frac{1}{12}\log\left|\dt\right| - \frac{1}{12}\log\left|\dt\right| + O(1) \\
           & = & \left(\frac{1}{3}+\frac{\delta}{2}\right)\log N + O(1) 
\end{eqnarray*}
as required.\\
}
\qed
\\

We now have to deal with the remaining case where $u \in \mathcal{P}_\tau$ has $\frac{1}{\sqrt{200}} \leq \left|u\right|$ which includes the cases of points on the components of $E(\R)$ not containing the identity (if they exist). 

\begin{lem}
\label{th:lem7}
Let $E$ be a curve given by $E: y^2 = x^3 + Ax + B$. Let $\left|\Delta_E\right|=N^\eta$ and suppose that $\left|A\right| \leq N^{4/3+2\delta}$ for some $\delta>0$. Suppose further that $E$ is associated with $\tau \in \c$, then for $N$ larger than an absolute constant,

\[ \text{Im}(\tau) \leq \max\left\{\frac{1}{2\pi}\left(4+6\delta - \eta\right)\log N +O(1) , D\right\} \]
\\
where $D > 0$ and $O(1)$ are absolute. 
\end{lem}
\proof{First suppose that $\tau \in \c_2$. Then Im$(\tau) \leq 1$ so this region is taken care of so long as we have $D \geq 1$. Now assume that $\tau \in \c_1 \cup \c_3$. We have

\[ \left|j(E)\right| = 1728\frac{\left|4A\right|^3}{\left|\Delta_E\right|} \leq CN^{4+6\delta-\eta} \]
\\
where  $C = 1728\times 64 $. Since $j(E) = j(E_\tau)$ we can use the $q$-expansion of $j(\tau)$

\[ j(\tau) = \frac{1}{q} + \sum_{n\geq0}c(n)q^n \]
\\
where $c(n) \in \Z_{>0}$. By a result of Petersson we have

\[ c(n) \sim \frac{e^{4\pi\sqrt{n}}}{\sqrt{2}n^{3/4}} \]
\\
as $n \rightarrow \infty$. Therefore there is some absolute $D > 0$ such that if $\text{Im}(\tau)=b \geq D$ then we have

\[ \left|j(\tau)\right| > \frac{1}{2\left|q\right|} = \frac{1}{2}e^{2\pi b} .\]
\\
Thus we have

\[ e^{2\pi b} < 2\left|j(E)\right| \leq(2C)N^{4+6\delta -\eta} .\]
\\
Allowing $N$ to be suitably large we arrive at our desired bound

\[ b < \left(\frac{4+6\delta-\eta}{2\pi}\right)\log N + \frac{\log 2C}{2\pi} = \left(\frac{4+6\delta-\eta}{2\pi}\right)\log N + O(1) .\]

}

\qed\\ 
\\
This lemma puts us in a position to deal with the remaining points.

\begin{lem}
\label{th:lem8}
Let $E$ be given by $E: y^2 = x^3 + Ax + B$ and suppose that $\left|A\right|\leq N^{4/3+2\delta}$ and $\left|\Delta_E\right| = N^{\eta}$ for some $\eta > 0$. If $E$ is associated with $\tau$ and $u \in \mathcal{P}_\tau$ is such that $\frac{1}{\sqrt{200}} \leq \left|u\right|$ then either
\[ \lambda_\infty(u) \leq \left(\frac{1}{3}+\frac{\delta}{2} - \frac{\eta}{12} \right)\log N + O(1)\]
\\
where the $O(1)$ is absolute or $\lambda_\infty(u)$ is bounded absolutely.
\end{lem}
\proof{We shall work by using the bound from Lemma ~\ref{th:lem3}:

\[ \lambda_\infty(u) \leq -\log\left|1-t\right| -\frac{1}{12}\log\left|\dt\right| + O(1) .\]
\\
First of all we wish to bound the term $-\log \left|1-t\right|$ absolutely. Let $u = a + b i$ so that we have 

\[ -\log\left|1-t\right| = -\log\left|1 - e^{2\pi i (a+bi)} \right| .\]
\\
Our aim is to show that $R=e^{2\pi i (a+bi)}$ is suitably far away from 1. Since $\left|u\right| = \left|a + bi\right| \geq \frac{1}{\sqrt{200}}$ we have that one of $\left|a\right|, \left|b\right| > \frac{1}{20}$. Suppose that $\left|a\right| \geq \frac{1}{20}$, then since the argument of the complex number $R$ is $2\pi a$ and $\left|a\right| \in [\frac{1}{20},\frac{1}{2}]$ we see that $\left|R - 1\right| \geq 2\sin\left(\frac{ \pi}{10}\right) > 0.3$. Thus in this case

\[ - \log\left|1-e^{2\pi i (a+bi)}\right| < -\log (0.3) < 1.21 .\]
\\
Suppose now that $\left|b\right| \geq \frac{1}{20}$. Then $|R| \leq e^{-\pi/10}$ and so $\left|R - 1\right| \geq 1 - e^{-\pi/10} > 0.27$. Therefore we have     

\[ - \log\left|1-e^{2\pi i (a+bi)}\right| < \log(0.27) < 1.31 .\]
\\
Therefore we see that $-\log\l|1-t\right| < 1.31$. Next, suppose $\tau \in \c_1 \cup \c_3$ and $\text{Im}(\tau)$ is greater than the $D$ in Lemma ~\ref{th:lem7}. Since Im$(\tau) > \frac{\sqrt{3}}{2}$ we see from the proof of Lemma ~\ref{th:lem3} that 

\[ \left|\dt\right| = \left|q\right| + O(1) \]
\\
where the $O(1)$ is absolute. So we have

\begin{eqnarray*}
-\frac{1}{12}\log\left|\dt\right| & = & -\frac{1}{12}\log\left|q\right| + O(1)\\
                                  & = & -\frac{1}{12}\log e^{-2\pi \text{Im}(\tau)} + O(1) \\
                                  & = & \frac{\pi \text{Im}(\tau)}{6} + O(1) \\ 
                                  & \leq & \frac{4+6\delta-\eta}{12}\log N + O(1) \\
                                  & = & \left(\frac{1}{3}+\frac{\delta}{2}-\frac{\eta}{12}\right)\log N +O(1) 
\end{eqnarray*}
where the penultimate inequality is from Lemma ~\ref{th:lem7}. This is the result we require. If $\text{Im}(\tau)$ is not larger than $D$ then it is uniformly bounded and thus so is $\left|\dt\right|$ which means that we have 

\[ \lambda_\infty(u) = O(1) .\]
\\
This completes the proof of the lemma.
}
\qed\\
\\
We are now in a position to prove Proposition ~\ref{th:thm3}.\\
\\
$\textit{Proof of Proposition ~\ref{th:thm3}:}$ Recall that we are trying to show that 

\[ \hat{h}(Q) \leq \left(\frac{1}{3} + \frac{\delta}{2} \right) \log N + O(1) \]
\\
where the $O(1)$ is absolute. We have the local height decomposition

\[ \hat{h}(Q) = \sum_{p \leq \infty}{\lambda_p (Q)}. \]
\\
Lemma ~\ref{th:lem1} tells us that 

\[ \sum_{p \neq \infty}{\lambda_p(Q)} \leq \frac{1}{12}\log \left|\Delta_E\right| .\]
\\
For the local height at the archemedian place, Lemma ~\ref{th:lem8} and Corollary ~\ref{th:cor2} both imply that we have either

\[ \lambda_\infty (Q) \leq \left( \frac{1}{3} + \frac{\delta}{2}\right) \log N - \frac{1}{12}\log \left|\Delta_E\right| + O(1) \]
\\
or we have 

\[ \lambda_\infty (Q) = O(1) \]
\\
where the $O(1)$ in both cases is absolute. In the first instance we can sum the local heights together to get the bound on $\hat{h}(Q)$ we desire. In the second instance we need to invoke the hypothesis that $\left|\Delta_E\right| \leq N^{4+6\delta}$ to attain the desired result. 
\qed

\section{The Large Sieve}

The Large Sieve deals with problems of the following kind. Let $S \subset [-N,N] \cap \Z$, $\p$ be a set of primes and $S_p$ be image of $S$ under the map $\Z \rightarrow \Z/p\Z$. Suppose that there is an $0 <\alpha < 1$ such that $\left|S_p\right|\leq \alpha p$, for every $p \in \p$. Then what can we say about $\left|S\right|$?\\
\\

\begin{lem}$\textbf{(The Large Sieve)}$ 
\label{th:lem9}
Let $S \subset \left\{M,\ldots,M+N\right\}$ for some $M \in \Z$ and $N>0$. Take $X>0$ and let $\p$ be a set of primes each no greater than $X$. Suppose that there is an $\alpha > 0$ such that $\left|S_p\right|\leq\alpha p$ for each $p \in \p$. Then 
\[ \left|S\right| \leq \frac{\alpha(N + X^2) }{(1-\alpha)\left|\p\right|} .\] 
\end{lem}

\proof{This follows from, e.g., \cite[Chapter~VIII,~Corollary~8.2.1, p.139]{CM}.}
\qed\\
\\
If we choose $N$ to be large enough so that if 

\[ \pi(x) = \#\left\{p \leq x : p \text{ is prime} \right\}  \]
\\
then

\[  \pi(N^{1/2}) > \frac{N^{1/2}}{2\log N^{1/2}} \]
\\
then by taking $\p$ in Lemma ~\ref{th:lem9} to consist of all primes between 43 and $N^{1/2}$ we see that 

\[ \left|\p\right| > \frac{N^{1/2}}{\log N} \]
\\
and so Lemma ~\ref{th:lem9} tells us that

\begin{eqnarray}\label{eq:5} \left|S\right| \leq \frac{\alpha 2 N^{1/2}\log N}{(1-\alpha)} .\end{eqnarray}
\\
The next lemma makes these endeavors relevant to our study of integral points on elliptic curves. 
	
\begin{lem}
\label{th:ls}
Let $E$ be an elliptic curve in Weierstrass form, $M \in \Z$ and let $S$ be the set

\begin{equation} S = \left\{ x(Q) \in \Z : \exists~Q=\left(x(Q),y(Q)\right) \in E(\Z) \right\} \cap [M,M+N] \nonumber\end{equation}
\\
for some $N >0$. Let $\eta > 0$ be such that $\left|\Delta_E\right| = N^{\eta}$. Then $\left|S\right| = O_{\eta}\left(N^{1/2}\log N\right)$. 
\end{lem}

\proof{
We are going to show this via the large sieve and in particular the inequality \eqref{e:3:15}. Let $p$ be a prime such that $p$ does not divide $\Delta_E$. Then $E$ is minimal over $\Q_p$ and the reduction of $E$ modulo $p$ is an elliptic curve $\tilde{E}_p$. By the Hasse Bounds (Theorem \ref{t:2:4}) we have that

\[ \left|\tilde{E}_p(\mathbb{F}_p)\right| \leq p + 1 + 2\sqrt{p}. \]
\\
If $p \not | 6$ then we can put $\tilde{E}_P$ into short Weierstrass form. With at most three exceptions, the points of $\tilde{E}_p(\mathbb{F}_p)$ pair up as a $\left\{Q, -Q \right\}$ with both points sharing the same $x$-coordinate. Thus

\[ \left|S_p\right| \leq \frac{p+1 + 2\sqrt{p}}{2} + \frac{3}{2} \]
\\
and this is less than $\frac{3}{4}p$ of $p > 42$. Let $\mathcal{P}$ be the set of primes less than $N^{1/2}$ excluding those $p$ less than 42 and those dividing $\Delta_E$. Since there are $O\left(\frac{\log N^{\eta}}{\log\log N^{\eta}}\right) = O_\eta(\log N)$ of these we see that we can apply inequality \eqref{e:3:15} with $\alpha = \frac{3}{4}$ to attain

\begin{equation}
\label{e:3:16}
\left|S\right| = O_{\eta}\left(N^{1/2}\log N\right).
\end{equation} \qed
}

The focus in this investigation is on uniformity and so we would like to remove the dependence on the size of $\Delta_E$ when dealing with points in $E(\Z,N)$. We can do this by appealing to a result due to Heath-Brown.

\begin{lem}
\label{l:3:3}
Let $F(\textbf{x}) \in \Z[\textbf{x}]$ be an irreducible form in three variables of degree $d$ with coprime coefficients. Let $\left\|F\right\|$ denote the maximum modulus of the coefficients of $F$. Let $\mathcal{C}: F = 0$. Then either $\C(\Q,N) \leq d^2$ or $\left\|F\right\| \ll N^{d(d+1)(d+2)/2}$.
\end{lem} 

\proof{
\cite[Theorem~4]{HB} \qed
}
\\
This result holds $\textit{a fortiori}$ for integral points on elliptic curves and shows that either $E(\Z,N)$ is bounded by $9$ or $\left\|E\right\| \ll N^{30}$. In the former case we have a bound far better than the Bombieri-Pila bound and so are happy. In the latter case we see that since $\left\|F\right\| \ll N^{30}$ we must have $\left|\Delta_E\right| \ll \left\|E\right\|^6 \ll N^{180}$. Thus we have an absolute bound on $\eta$ which means that for our set up we can make the implied constant in Lemma \ref{th:ls} absolute.\\

\begin{Prop}
\label{th:lsprop}
Let $E$ be an elliptic curve given by $y^2 = x^3 + Ax + B$ with $A,B \in \Z$. Let $0 < M \leq N$ and $\e>0$. Suppose that either $M>N^{2/3+2\e}$ or $\left|A\right| > N^{4/3+4\e}$. Then 

\[ \#\left(E(\Z) \cap (([-N,-M]\cup [M,N])\times[-N,N])\right) = O(N^{1/3-\e/2}). \] 
\end{Prop}

\proof{ 
We can use Lemma ~\ref{th:ls} to attain such a bound if we were to know that the portion of $E(\R)$ which lies over the intervals $[-N,-M]\cup [M,N]$ in fact lies over a small absolute number of intervals each of length less than $N^{2/3-\e}$. This phenomenon occurs for the portions of $E(\mathbb{R})$ which have gradient of absolute value exceeding $N^{1/3 + \e}$. Thus we should investigate the gradient of $E$. We have

\[ \frac{\text{dy}}{\text{dx}} = \frac{3x^2 + A}{2y}. \]
\\
We are interested in the absolute value of the gradient being larger than $N^{1/3+\e}$ which occurs when we have

\[ \left|3x^2 + A\right| \geq 2\left|y\right|N^{1/3+\e}. \]
\\
Since $\left|y\right| \leq N$ see that this inequality holds if we have

\begin{equation}
\label{eq:lsineq}
\left|3x^2 + A\right| \geq N^{4/3+\e} . 
\end{equation}
\\
This is the inequality which we will study. If $A \geq 0$ then ~\eqref{eq:lsineq} holds if either $A \geq N^{4/3+\e}$ or $\left|x\right| \geq N^{2/3+\e/2}$. The hypotheses of the proposition clearly imply that we are in one of these cases and so we are done for the case $A \geq 0$. \\
\\
Now suppose that $A < 0$. Let us begin by rewriting the equation for $E$ as

\[E: y^2 = x^3 - Cx + B \]
\\
where $C = \left|A\right|$ so that the sign of the symbols is clearer. Suppose that $C = \left|A\right| \leq N^{4/3 + 4\e}$. Then $M \geq N^{2/3+2\e}$ whence

\[ 3x^2 - C \geq 3N^{4/3+4\e} - C \geq 2 N^{4/3+4\e}. \]
\\
In this case we see that ~\eqref{eq:lsineq} holds and so we are done. Thus we can assume that $C > N^{4/3 + 4\e}$. Suppose that ~\eqref{eq:lsineq} fails. This would occur if $3x^2$ were too close to $C$ so we shall write $x = \pm \sqrt{\frac{C}{3}} + t$. Then we have

\begin{eqnarray}
\label{eq:ls2}
2N^{4/3+\e} & > & \left|3x^2 - C\right| \nonumber \\
            & = & \left|3 \left(\pm\sqrt{\frac{C}{3}}+t\right)^2 - C \right| \nonumber\\
            & = & \left|6t\left(\pm \sqrt{\frac{C}{3}} \right) + 3t^2 \right|. 
\end{eqnarray}  

Let $R = \pm \sqrt{\frac{C}{3}}$. If $Rt \geq 0$ then ~\eqref{eq:ls2} implies that

\[ 6\sqrt{\frac{C}{3}}\left|t\right| < 2N^{4/3+\e} .\]
\\
Recalling that $C > N^{4/3+4\e}$, we have

\[ \left|t\right| < \left|\frac{2\sqrt{3}N^{4/3+\e}}{6\sqrt{C}}\right| < N^{2/3-\e} .\]
\\
This means that $x$ lies in a interval of length at most $2N^{2/3-\e}$ which allows us to apply Lemma ~\ref{th:ls} to attain that there are no more than $O\left(N^{1/3-\e/2}\right)$ such integer points. Let us now suppose that $Rt < 0$. We split into two cases: $x \geq 0 $ and $x < 0$. In the first case we can write $x = \sqrt{\frac{C}{3}} + t$ where $-\sqrt{\frac{C}{3}} \leq t \leq 0$. Hence

\[ \left| 6 \sqrt{\frac{C}{3}}t + 3t^2\right| \geq 3\sqrt{\frac{C}{3}}\left|t\right| \]
\\
and putting this into \eqref{eq:ls2} above yields

\[ \left|t\right| \leq \frac{2N^{4/3+\e}}{\sqrt{3C}} < \frac{2}{\sqrt{3}}N^{2/3-\e}. \]
\\
As before, this fact combined with Lemma \ref{th:ls} tells us that there are at most $O\left(N^{1/3-\e/2}\right)$ such integer points. This just leaves the case of $x < 0$. In this case we have $x = -\sqrt{\frac{C}{3}} + t$ where $0 \leq t \leq \sqrt{\frac{C}{3}}$. Thus

\[ \left| -6 \sqrt{\frac{C}{3}}t + 3t^2\right| \geq 3\sqrt{\frac{C}{3}}\left|t\right| \]
\\
and so, as above, we have

\[ \left|t\right| \leq \frac{2}{\sqrt{3}}N^{2/3-\e} .\]
\\
Hence we see that in this case as well there are no more than $O\left(N^{1/3-\e/2}\right)$ such integer points. We are now done since the portion of $E(\mathbb{R})$ in the range given in the statement of the proposition can be split up into 3 connected arcs where the gradient is greater than $N^{1/3+\e}$ (3 arcs since the degree of $E$ is 3), these each lie over a region on the $x$-axis of length no more than $2N^{2/3-\e}$ and the remaining part of the curve can be split into 3 regions corresponding to the cases which we dealt with above (i.e, the sign of $A$ and the sign of $Rt$). Each of these regions lies above an interval of length less than $ \frac{4}{\sqrt{3}}N^{2/3-\e}$. Applying Lemma ~\ref{th:ls} to each of these regions and summing them together gives us the desired result.
\qed 
} 
\section{Proof of Main Theorem}

In this section we shall prove the Main Theorem first by showing it holds for curves in short Weierstrass form and then by reducing the more general case to this form. First of all, let $\delta$ be the same as in Theorem ~\ref{th:thm2}. If we look at points $Q \in E(\Z)\cap [-N,N]^2$ such that $\left|x(Q)\right|>N^{2/3+2\delta/k}$ for some $k>2$ then Proposition \ref{th:lsprop} taken with $\e = \delta/k$ tells us that there are no more than

\[ O\left(N^{1/3-\delta/2k}\right) \]
\\
of these. Furthermore, if $\left|A\right|>N^{4/3+4\delta/k}$ then we can again apply Proposition ~\ref{th:lsprop} with the same $\e$ to get the same bound.  \\
\\
Now, suppose that $\left|A\right|\leq N^{4/3+4\delta/k}$ and we restrict to points $Q \in E(\Z)\cap [-N,N]^2$ with $\left|x(Q)\right|< N^{2/3+2\delta/k}$ then we see that

\[ \left|B\right| \leq \left|y^2\right| +\left|x^3\right| + \left|Ax\right| \leq N^{2+ 6\delta/k} \]
\\
and so $\left|\Delta_E\right| \leq N^{4+6(2\delta/k)}$. This allows us to apply Proposition ~\ref{th:thm1} with the $\delta$ of that Proposition taken to be $2\delta/k$ which in combination with the bounds of Helfgott and Venkatesh (Theorem ~\ref{th:thm2}) tell us that the number of integer points in this case is bounded by

\[ O_\e(N^{(1/3+\delta/k)(1-\delta)+\e}) = O_\e(N^{1/3-\delta(k-3)/3k - \delta^2/k}) .\]
\\
Since we have $k>2$ we see that the above exponent is less than $\frac{1}{3}$. We can therefore take the maximum of $\frac{1}{3} - \frac{\delta}{2k}$ and $\frac{1}{3}-\delta \frac{k-3}{3k} - \frac{\delta^2}{k}$ as the exponent for a bound on $E(\Z,N)$. Since we are not concerned with optimisation here we shall be content with noting that if we ignore the $\delta^2/k$ term (which only acts in our favour) then the best value for $k$ is $\frac{9}{2}$ which gives us that $E(\Z,N) \ll O(N^{1/3 - \delta/9})$. In any case we have the following

\begin{Prop}
\label{th:thm5}
Let $E$ be an elliptic curve given by the equation $y^2 = x^3 + Ax + B$ where $A,B \in \Z$. Then there is a $\delta >0$ such that

\[ \# E(\Z,N) = O(N^{1/3-\delta}) .\]
\end{Prop}
In order to prove the Main Theorem we shall need a result drawing on the work of Ellenberg and Venkatesh in \cite{EV}
	
\begin{lem}
\label{th:lem11}
Let $f \in Z[X,Y]$ be a homogeneous polynomial of degree $d$, let $Z \subset \mathbb{P}^2$ be the curve cut out by $f$. Then the number of rational points on $Z$ of naive height no more than $N$ is bounded by 
	
\[ (N^{2/d}\left\|f \right\|^{-1/{d^2}}+1)N^{\e} \] 
\\
where $\left\|f\right\|$ denotes the $\ell^2$ norm.   
\end{lem}  
	
\proof{This is the $n=1$ case of \cite[Proposition~2.1]{EV} }
\qed \\
	
\begin{lem}
\label{th:lem12}
Let $E,f$ be as in the statement of the Main Theorem so $E: y^2 = x^3 + Cx^2 + Dx + F = f(x)$. Then if $\left|C\right|\geq N^{1+6\e}$ we have that 
	
\[ \# E(\Z,N) = O\left(N^{1/3 - \e}\right) .\]
\end{lem}
\proof{We shall use Lemma ~\ref{th:lem11} in a fairly crude way. Let $g(x,y) = y^2 - f(x)$. To $g$ we associate the homogeneous form
	
\[ G(X,Y,Z) = Y^2 Z - X^3 - CN^2 X^2 Z - D N X Z - F Z^3 .\]
\\
(This homogenisation is found in \cite[Remark~2.2]{EV}). Clearly the $\Z$ points of $E$ in our box correspond to points on $G(X,Y,N) = 0$ of height no greater than $N$. Hence we can use the bounds of Lemma ~\ref{th:lem11} to get the bound 
	
\begin{eqnarray}\label{eq:9}  \# E(\Z,N) \leq (N^{2/3}\left\|G \right\|^{-1/{9}}+1)N^{\e} .\end{eqnarray}
\\
Now, if $\left|C\right| \geq N^{1+6\e}$ we see that $\left\|G\right\| \geq N^{3+6\e}$ and so ~\eqref{eq:9} becomes
	
\[  \# E(\Z,N) \leq (N^{2/3}N^{-1/{3}-2\e}+1)N^{\e} = O\left(N^{1/3-\e}\right) \]
\\
as required.
\qed}\\
\\
$\textit{Proof of Main Theorem.}$ We begin by changing the coordinates in order to bring the curve first into a short Weierstrass form. We move from coordinates $(x,y)$ to $(X,Y)$ where 

\[ X = 9x + 3C ~~~~~~~~~~~ Y = 27y .\]
\\
This gives us a new expression for $E$

\[ E^{\prime} : Y^2 = X^3 + Ax + B \]
\\
with $A,B \in \Z$. Clearly integral points of $E$ are sent to integral points on $E^\prime$. Moreover, points in $\# E(\Z,N)$ are sent to points of $E$ in the box 

\[ \mathcal{B} = [-9N + 3C, 9N + 3C] \times [-27N,27N].\]
\\
Now, Lemma ~\ref{th:lem12} allows us to assume that $\left|C\right| \leq N^{1+6\e}$ for some $\e$ for otherwise we have the bound
	
\[ O\left(N^{1/3 - \e} \right) \] 
\\	
and so could take $\delta^\prime = \e$.  Thus we have that 
	
\[\mathcal{B} \subset [-4N^{1+6\e},4N^{1+6\e}]^2 .\] 
\\
By Proposition ~\ref{th:thm5} we know that there is a $\delta > 0$ such that the number of integral points in this box is 
	
\[ O\left( N^{(1+6\e)/3 - \delta} \right) \]   
\\
and so by choosing $\e$ small enough we see that there is a $\delta^\prime > 0$ such that

\[ \# E(\Z,N) = O\left(N^{1/3 - \delta^\prime}\right) .\qed\]

\section{Extending the result}

So far we have been considering the number of integer points on an elliptic curve which fall inside a square box centred at the origin. Our Main Theorem applies in the case that the curve is of the shape $y^2 = f(x)$. In this section we shall consider the problem of placing the box arbitrarily in the plane for curves of this shape.\\
\\
As before, we shall start by considering curves in short Weierstrass form. Let $E$ be an elliptic curve in short Weierstrass form, $N>0$ and $Q = (x(Q),y(Q)) \in E(\Z)$. Let $\alpha, \eta >0$ be such that $\left|x(Q)\r| = N^\alpha$ and $\l|y(Q)\r| = N^\eta$. The main observation we shall exploit is contained in the following

\begin{lem}
\label{th:lem13}
Let $E$ be an elliptic curve in short Weierstrass form with $\left|j(E)\right| > \e_1$ and $\left| A\right| \leq N^{4\eta/3 - \e_0}$ for some $\e_0,\e_1 > 0$. Let $Q, \al, \eta$ be defined as above. Then for any $\e > 0$ and $N>0$ sufficiently large with respect to $\e, \e_0,\e_1$ we have

\[ \frac{\eta}{3} \leq \frac{\alpha}{2} + \e .\] 
\end{lem}  

\proof{
Since $E$ has the form

\[ y^2 = x^3 + Ax + B\]
\\
we easily see that

\[N^{2\eta} \leq N^{3\al} + \left|A\right|N^{\al} + \left|B\right| \leq N^{3\al} + N^{\al + 4\eta/3 - \e_0} + \l| B\r|.\]
\\
Suppose for a moment that we can assume that $\l| B \r| \leq D N^{2\eta - \kappa}$ for some $D,\kappa >0$. Then for $N$ sufficiently large we would have $N^{-\kappa} < \frac{1}{2D}$ and so 

\[ N^{2\eta} \leq 2\left(N^{3\al} + N^{\al + 4\eta/3 - \e_0}\right) .\]
\\
It is not hard to see that the lemma would follow from such an inequality so long as $\kappa$ and $D$ are bounded in terms of $\e,\e_0,\e_1$. Thus we just need to control $\l|B \r|$. Since $\left|j(E)\right| > \e_1$ we have that

\begin{equation}\label{eq:b} 1728(4\left|A\right|^3) \geq \e_1\left|4A^3 + 27B^2\right| .\end{equation}
\\
If $A \geq 0$ then ~\eqref{eq:b} yields

\[ \left|B\right| \leq 2\sqrt{\frac{1728-\e_1}{27\e_1}}A^{3/2} \leq D_{\e_1} N^{2\eta - 3\e_0 /2} \]
\\
which is as we wanted. Now suppose that $A < 0$. Set $C = \left|A\right|$. Then from ~\eqref{eq:b} we have

\[ 1728(4C^3) \geq \e_1\left|27B^2 - 4C^3\right| .\]
\\
Suppose that $27B^2 \geq 4C^3$. Then we have

\[ \left|B\right| \leq 2\sqrt{\frac{1728+\e_1}{27\e_1}}C^{3/2} \leq D_{\e_1} N^{2\eta - 3\e_0 /2} \] 
\\
as required. If, on the other hand, $27B^2 \leq 4C^3$ then

\[ \left|B\right| \leq \frac{2C^{3/2}}{3\sqrt{3}} \leq \frac{2}{3\sqrt{3}}N^{2 -3\e_0 /2} \]
\\
which is as we want.
\qed\\}
\\
We shall now use this lemma to prove Theorem ~\ref{th:thm1}.\\

\newtheorem{Thmm}{Theorem}
\begin{Thmm}
Let $E$ be an elliptic curve given by the equation $y^2 = f(x)$ where $f \in \Z[x]$ is a monic cubic with no repeated roots. Suppose further $c_4(E) < 0$ and $j(E) > \e$ for some $\e > 0$. Let $\mathcal{B}$ be any box in the plane with sides of length $N$. Then there is a $\delta > 0$ (depending only on $\e$) such that $\# \left(E(\Z) \cap \mathcal{B}\right) = O(N^{1/3-\delta})$.  
\end{Thmm}  
\proof{

We shall begin by assuming that $E$ is in short Weierstrass form and then derive the more general statement from this. So, let $E$ be cut out by the equation

\[ E: y^2 = x^3 + Ax + B\]
\\
and let $\mathcal{B}$ be a box of size $N$ by $N$ somewhere in the plane. The aim is to get a bound on $\#E(\Z) \cap \mathcal{B}$ of the quality $O(N^{1/3-\e})$ for some $\e > 0$. We shall be using the assumption that $A > 0$ (which is the same as assuming $c_4(E) < 0$). We would like to consider where gradient considerations and the Large Sieve can give us such a result via Lemma \ref{th:ls}. To do this we need to be able to control the size of $\Delta_E$. Suppose that $\mathcal{B}$ is centred at the integer point $(x_0,y_0)$. Then by the change of variables

\[ X = x - x_0 ~~~~~~~~~Y = y - y_0 \]
\\
we can move $\mathcal{B}$ to be centred at the origin. This has the effect of changing our equation to

\[ E^\prime : Y^2 + 2Yy_0 + y_0^2 = X^3 + 3X^2x_0 3Xx_0^2 + x_0^3 + AX + Ax_0 + B .\]
\\
It is clear then that $\left\|E^\prime\right\| \geq \max\{ \left|y_0\right| , \left|x_0\right \}$ and so by Lemma \ref{th:lem11} we see that if either $\left|x_0\right|$ or $\left|y_0\right|$ are greater than $N^{3+18\e}$ then $E^\prime(\Z,N) \ll N^{1/3-\e}$ which implies the result of this Theorem. Hence we may assume that both are bounded by $N^{3+18\e}$. In this case we have that the points $E(\Z)\cap \mathcal{B}$ (note that this is $E$ not $E^\prime$) are counted by $E(\Z,2N^{3+18\e})$. Applying Heath-Brown's result (Lemma \ref{l:3:3}) we see that we have $\left\|E\right\| \ll N^{30(3+18\e)}$ which in turn means that $\left|\Delta_E\right| \ll N^{180(3+18\e)}$. This bound for the size of $\Delta_E$ is good enough for our needs since the exponent of $N$ is absolute. We can now proceed with the use of the Large Sieve technique. We have

\[ \left|\frac{\text{dy}}{\text{dx}}\right| = \left|\frac{3x^2 + A }{2y}\right| \geq \left|\frac{3x^2 + A}{2N^{\eta}}\right| .\]
\\
Since our box is of size $N$ by $N$ we see by the same reasoning as before that the Large Sieve delivers the bound we want if we have $\left|\frac{\text{dy}}{\text{dx}}\right| \geq N^{1/3 + 2\e}$ hence we shall assume that this is not the case. Thus we have

\[ \left|3x^2 + A\right| \leq 2N^{1/3+\eta+\e}. \]
\\
Since $A > 0$ we learn that

\[ \left|x\right| \leq N^{1/6 + \eta/2 + \e/2} \]
\\
and

\[ \left|A\right| \leq 2N^{1/3 + \eta + \e} .\]
\\
We wish to apply Lemma ~\ref{th:lem13} and so need to check that the hypothesis on $A$ is satisfied. This would be the case if

\[ \frac{1}{3} + \eta + \e \leq \frac{4\eta}{3} - \e_0 \]
\\
for some $\e_0 > 0$. This is true if and only if

\[ \eta \geq 1 + 3(\e + \e_0) .\]
\\
If this inequality were false then by a change of coordinates we can move the box horizontally so that it lies over the interval $[-N,N]$. This would have the effect of putting the curve into the form $y^2 = f(x)$ and so we could apply the Main Theorem with a box of size $N^{1+3(\e+\e_0)}$. So long as $\e,\e_0$ were small enough (with respect to the $\delta$ in the Main Theorem) we could get the type of bound we want. Thus we see that we may assume the hypothesis on $A$ required for Lemma ~\ref{th:lem13} safely. We thus apply the lemma with $\al = \frac{1 + 3\eta + 3\e}{6}$ to attain

\[ \frac{\eta}{3} \leq \frac{1}{12} + \frac{\eta}{4} + \frac{5\e}{4} \]
\\
and so

\[ \eta \leq 1 + 15\e \]
\\
which, by the same argument as above, can be dealt with by the Main Theorem for $\e$ sufficiently small. This concludes the proof of the Theorem for the case where $E$ is in short Weierstrass form.\\
\\
Suppose now that $E$ is given by $y^2 = f(x)$ where $f$ is as in the statement of the Theorem. Let $\mathcal{B}$ be a box of size $N$ by $N$ anywhere in the plane. By an obvious change of coordinates we can move $\mathcal{B}$ so that it is centred at the origin. This has the effect of changing the equation for $E$ into something of the form 

\[ y^2 + ay = g(x) \]
\\
where $a \in \Z$ and $g \in \Z[x]$ is a monic cubic. By the same application of the results of Ellenberg and Venkatesh (Lemma \ref{th:lem11}) as we used in the proof of the Main Theorem we can replace $g$ with a cubic of the form $x^3 + Ax + B$ and at the expense of increasing the width of $\mathcal{B}$ to $4N^{1+6\e}$ for some $\e>0$ which we can make as small as we like. By applying one last change of variables we can eliminate the $ay$ term which results in the box being moved vertically and replaces our equation with one in short Weierstrass form. This has been dealt with already and so we are done. 
}\qed\\
\\
We shall end this section by making a few remarks about why Theorem \ref{th:thm1} has the hypotheses it does and why we cannot remove them at present. Firstly, why the shape of the equation is $y^2 = f(x)$ rather than the more general Weierstrass form of $y^2 + a_1xy + a_3y = f(x)$. It is clear that we can deal with equations of the form $y^2 + a_3y = f(x)$ with Theorem \ref{th:thm1} (so long as the other hypotheses are met) as was seen above in the proof of Theorem \ref{th:thm1}. The reason we cannot deduce results on the long form from the form in Theorem \ref{th:thm1} is that the change of variables 

\[ Y = y + \frac{a_1}{2}x ~~~~,~~~~X=x \]
\\
which carries the equation $y^2 + a_1xy + a_3y = f(x)$ to one of the form $Y^2 + a_3Y = F(X)$ has the effect of changing the box $[-N,N]^2$ (in which we are performing our count) to a parallelogram which has width along the $X$-axis of $2N$ but height on the $Y$-axis of $O(a_1N)$. We can get a control on $\left|a_1\right|$ from the bounds of Ellenberg and Venkatesh as we did in Lemma \ref{th:lem11} yielding that we have an improvement on the Bombieri-Pila bounds if $\left|a_1\right| > N^{1+\e}$ for some $\e >0$. However, this control is not strong enough to get a result in the case $\left|a_1\right| < N^{1+\e}$ by change of variables. This is because the parallelogram mentioned above could have vertical height up to $O(N^{2+\e})$ which is too large for our methods to apply.\\
\\
As for the hypotheses $\left|j(E)\right| \geq \e$ and $c_4(E) < 0$, these can both be understood in terms of our application of the large sieve. The idea is that the large sieve helps when the gradient of $E(\R)$ is particularly steep. Let us consider the case where $j(E) = 0$. These are the Mordell curves $E_D: y^2 = x^3 + D$. If we consider the cases where $D > 0$ and indeed is substatially larger than $N$ then when we look at the portion of $E_D(\R)$ which lies above the interval $[-N,N]$ on the $x$-axis, we see that it has a rather flat slope. It turns out that for many values of $D$ (large in terms of $N$) this slope is too shallow for us to apply the large sieve and so we cannot get our result. The same problem arises for curves with $c_4(E) < 0$. It may be that there is a more subtle way of using the sieve which can help here but it is my opinion that a new method is needed to deal with both these hypotheses and the case of the general Weierstrass form.

\section{An application to degree 1 del Pezzo surfaces}

A del Pezzo surface of degree $d$ defined over $\Q$ is a smooth, projective algebraic surface which has ample anticanonical bundle $-K_X$ and where the self-intersection of $K_X$ with itself has intersection number $K_X^2 = d$. For these surfaces we may take the height $H_X$ associated with the anticanonical embedding and study the counting function

\[ X(\Q,N) = \# \left\{ P \in X(\Q) : H_X(P) \leq N \right\} .\]
\\
The main conjecture in this area is Manin's Conjecture \cite{BM} which in the case of del Pezzo surfaces claims that there is a Zariski open subset $U \subset X$ so that one has

\[ U(\Q,N) \sim c_X N(\log N)^{\rho_X -1} \]
\\
where $c_X$ is a constant depending on $X$ and $\rho_X$ is the Picard rank of $X$. There is an industry in proving Manin's Conjecture for specific del Pezzo surfaces, often using delicate arguments from analytic number theory (for an overview see \cite{Br}) and from experience there it is understood that the problem becomes more difficult the lower the degree of the del Pezzo surface. Evidence for this comes from the fact that Manin's Conjecture has not been verified for a single degree 1 del Pezzo surface. In terms of upper bounds for $X(\Q,N)$, little is known in general. That the results proved thus far in this paper can serve to do this was pointed out to me independently by Daniel Loughran and Pierre Le Boudec.\\
\\
\begin{Prop}
Let $X$ be a degree 1 del Pezzo suface defined over $\Q$ and let $N>0$. Then there is a $\delta > 0$ such that $X(\Q,N) \ll N^{3 - \delta}$ where the implied constant is uniform in $X$.
\end{Prop}

\proof{
It is known (eg, \cite[\textsection~2.2.3]{TA}) that $X$ can be embedded into weighted projective space $\mathbb{P}_{\Q}(3,2,1,1)$ as a smooth sextic of the shape 

\begin{equation}
\label{e:dp1}S: y^2 = x^3 + F_4(u,v)x + F_6(u,v)
\end{equation}
\\
where $y,x,u,v$ have weights $3,2,1,1$ respectively and where $F_i$ are homogeneous polynomials of degree $i$. In this setting, the counting function $X(\Q,N)$ is bounded by the size of the set

\[ S(N) = \left\{ P= (y,x,u,v) \in S(\Q) : x,y,u,v \in \Z, \left|x\right| \leq N^2, \left|y\right| \leq N^3, \left|u\right| \leq N, \left|v\right| \leq N \right\}. \]
\\
For each choice of a pair $(u,v)\in \Z^2$ with $\left|u\right|,\left|v\right| \leq N$ (of which there are clearly $O\left(N^2\right)$) we have that \eqref{e:dp1} specialises to the equation for a cubic curve in short Weierstrass form with integer coefficients which we shall denote $E_{(u,v)}$. This curve will be an elliptic curve so long as we have 

\[ \Delta(E_{(u,v)}) = -16(4F_4(u,v)^3 + 27F_6(u,v)^2) \neq 0 .\]
\\
The pairs for which $\Delta(E_{(u,v)}) = 0$ are integer points lying on a degree 12 curve in $\mathbb{A}^2_{\Q}$ and so by the Bombieri-Pila bounds the number of these points for which $(u,v) \in [-N,N]^2$ is bounded by $O_\e(N^{1/12+\e})$. Let $(u,v)$ be such a pair and consider then the curve $E_{(u,v)}$. All of the points $(x,y) \in E_{(u,v)}(\Z)$ for which $(y,x,u,v)$ lie in $S(N)$ are within the box $[-N^3,N^3]^2$ and so applying Bombieri-Pila again we see that $E_{(u,v)}(\Z,N^3) \ll_\e N^{1 + \e}$. Thus these curves contribute at most $O_\e\left( N^{13/12 + \e} \right)$ points to $S(N)$. Since we are interested in a bound of the quality $O\left(N^{2-\delta}\right)$ then we see that these points are insignificant to the count so long as $\delta < 11/12 - \e$ (which it will be).\\ 
\\
Now, let $(u,v)$ be a pair such that $E_{(u,v)}$ is an elliptic curve. By the Main Theorem we see that $E_{(u,v)}(\Z,N^3) \ll N^{1-\delta}$ for some $\delta > 0$. Taking into account the different choices for the pair $(u,v)$ we attain the desired bound on $X(\Q,N)$. Of course this $\delta$ is thrice that of the $\delta$ in the Main Theorem. \qed\\
}
\\
In certain special cases then there are some results due to Munshi in \cite{RM} who shows that for del Pezzo surfaces of the form

\begin{eqnarray*}
S_1 : & y^2 = x^3 + F_4(u,v)x \\
S_2 : & y^2 = x^3 + b(F_2(u,v))^3 \\
S_3 : & y^2 = x^3 + b(F_3(u,v))^2 
\end{eqnarray*}
\\
where $b\in \Z$ then one can obtain $S_i(\Q,N) \ll_b N^{2+\e}$ where the implied constant depends on $b$ for the last two examples but not on the $F_i$ in any case. These follow from the fact that we have a non-trivial torsion section of these elliptic surfaces over either $\Q$ in the first case, $\Q(\sqrt[3]{b})$ in the second and over $\Q(\sqrt{b})$ in the third. Such an approach will not work in general since we cannot hope for a torsion section of $X$ to exist over a fixed number field $K/\Q$ in general. We also mention that through a different approach Munshi also shows that for 

\[ X : y^2 = x^3 + (F_2(u,v))^3 \]
\\
one can obtain $X(\Q,N) \ll N^{4/3+\e}$. We note that whilst our results are not as strong, they are hold for any such surface and are completely uniform in terms of the surface $X$. \\

\section{Concluding Remarks}

It would be desirable to remove the hypotheses in Theorem ~\ref{th:thm1} and to extend it curves in long Weierstrass form however the methods employed in this paper do not seem to be able to handle these extensions. We also note that our methods make use of a lot of information regarding elliptic curves, most importantly the explicit equations for the local canonical heights. This means that the generalisation of the methods seen here to higher genus curves would not be trivial. It has been suggested that the correct setting for such arguments would be within Arakelov Theory. It would be interesting to see if such a reformulation of these ideas in that domain could lead to generalisations of the results seen here both to number fields and to higher genus curves.\\

\end{document}